\documentclass{pcmAX}
\usepackage{graphicx} 
\usepackage{float}
\usepackage{polski}
\usepackage[polish,english]{babel}
\usepackage{bbm}
\usepackage{amsmath}

\usepackage[latin2]{inputenc}

\author{Maciej Wielgus}

\nralbumu{262778}

\title{Perona-Malik equation and its numerical
properties}

\tytulang{Perona-Malik equation and its numerical
properties}

\kierunek{Matematyka}

\opiekun{dra Pawła Bechlera\\
  Instytut Matematyki Stosowanej i Mechaniki\\
  }

\date{Lipiec 2010}

\begin{document}
\maketitle

\begin{abstract}
 	
This work concerns the Perona-Malik equation, which plays essential role in image processing. The first part gives a survey of results on existance, uniqueness and stability of solutions, the second part introduces discretisations of equation and deals with an analysis of discrete problem. In the last part I present some numerical results, in particular with algorithms applied to real images. 
\end{abstract}

\renewcommand{\contentsname}{Table of contents}
\renewcommand{\bibname}{Bibliography}
\tableofcontents


\chapter*{Introduction}
\addcontentsline{toc}{chapter}{Introduction}

\ \ \ The Perona-Malik model, first proposed in 1987  \cite{PM}, is a nonlinear partial diffusion equation that uses an inhomogeneous diffusivity coefficient. It is widely used in image processing for purposes like smoothing,  restoration, segmentation, filtering or detecting edges. 
\\
\\

It was shown by Kawohl and Kutev that the equation may have no global weak solutions in $C^1$ \cite{KK}. There are several methods which deal with this problem by regularizing the equation, both in space and time (see for example \cite{CLMC}). 
\\
\\

Although the basic model is ill-posed, its discretizations are found to be stable. This fact is sometimes referred to as the Perona-Malik paradox \cite{Kit}. The explanation for these observations was given by Weickert and Benhamouda [WB], who investigated the regularizing effect of a~standard finite difference discretization.
\\
\\

The following work presents an analysis of the Perona-Malik equation, surveys the results on ill-posedness of the continuous problem and its discrete versions. I also present numerical results, compare different versions of the algorithm and demonstrate the evolution of a picture under the Perona-Malik filter. As I present the Perona-Malik model from the perspective of image processing, I usually bound myself to the two-dimensional version of the equation, or, for clarity of demonstration, show the results in one-dimension. I do not include proofs as they are to be found in the publications listed in the bibliography.
\\
\\

All implementations are done in the Octave/Matlab environment.

\chapter{The Perona-Malik equation}\label{r:PME}


\section{Gaussian diffusion in image processing}

The diffusion equation is a partial differential equation often arising from physical phenomena. It is usually presented in a general form

\begin{equation}
 \partial_{t} u(x,t) = \nabla( C(u,x,t) \nabla u(x,t))
\end{equation}

where the coefficient $C$ is called \textit{diffusivity} and, depending on the model, can be a scalar, scalar function of coordinates (inhomogeneity) or a tensor (anisotropy). Moreover, the equation becomes nonlinear if the diffusivity coefficient depends on the solution $u$.
\\

For a constant $C \equiv 1$ we obtain a simple linear model of the heat equation, which leads to the following Cauchy problem

\begin{equation}
\label{Heateq}
 \left\{
\begin{array}{ll} 
u_t = \Delta u & \mathrm{for} \ x \in \mathbbm{R}^n, t>0 \\
u(x,0) = f(x) & \mathrm{ for} \ x \in \mathbbm{R}^n \\

\end{array} \right.
\end{equation}
\\

\textbf{Theorem 1.1} (\cite{Strzel})  \textit{If $f \in C^0( \mathbbm{R}^n ) \cap L^{\infty}( \mathbbm{R}^n)$ then}

\begin{equation}
\label{Heatsol}
u(x,t) = \left\{
\begin{array}{ll} 
\int_{\mathbbm{R}^n} f(y)G(x-y,t) \ dy = f * G_t & \mathrm{for} \ t>0 \\
f(x) & \mathrm{ for} \ t=0 \\

\end{array} \right.
\end{equation}

\textit{where G is Gaussian function}
\begin{equation}
 G(x,t) = \frac{\exp(-|x|^2/4t) }{2^n (\pi t)^{n/2} }, \quad x \in \mathbbm{R}^n , \ t>0, \end{equation}

\textit{is the solution to \eqref{Heateq}. It is unique under the assumption that}

\[ \forall T > 0 \ \ \exists \  C_T, a_T > 0 : \ \ |u(x,t)| \le C_T \exp(a_T |x|^2)\quad \mathrm{for} \ \ (x,t) \in \mathbbm{R}^n \times [0,T]  \]

\textit{or} 

\[ u(x,t) \geq 0 \quad \mathrm{for} \ \ (x,t) \in \mathbbm{R}^n \times [0,T] \]

The second condition for uniqueness is particularly interesting, because of its physical interpretation  (temperature is never below absolute zero).
\\

It is visible from \eqref{Heatsol}, that solving the heat equation means convolving the initial temperature distribution with the Gauss function with the standard deviation parameter $\sigma = \sqrt{2t}$. Such operation is well known as an efficient down-pass spectral filtration and it is commonly used for smoothing pictures by averaging values within a certain neighbourhood. This is where the basic idea of using diffusion equations in image processing comes from.
\\

An important property of the heat equation arises from the properties of convolution:

\begin{equation}
 \int_{\mathbbm{R}^n} u(x,t) \ dx =  \int_{\mathbbm{R}^n}f * G_t \ dx = \int_{\mathbbm{R}^n} f \ dx \cdot \int_{\mathbbm{R}^n} G_t \ dx = \int_{\mathbbm{R}^n} f \ dx \cdot 1 = \int_{\mathbbm{R}^n} u(x,0) \ dx
\end{equation}

The interpretation of this result would be the energy conservation law for the transport of heat, or, in image processing, the conservation of the average intensity of the given image.
\\

\section{Introducing the Perona-Malik equation}

Using  Gaussian diffusion for image filtration has one great drawback: its homogeneity, which leads to unwanted diffusion of image features like edges.
The idea behind the Perona-Malik equation is to modify the heat equation by adding the diffusivity coefficient depending on space activity in a given part of a picture, measured by the norm of the local picture gradient. For small gradient norm (homogeneous regions) large values of the diffusivity are expected, to perform stronger smoothing. In regions with big gradient norm (inhomogeneity) smaller diffusivity is expected, to slow down the diffusion process and protect delicate image features. Hence, the precise Perona-Malik problem formulation with a Neumann boundary condition is

\begin{equation}
\label{PMeq}
 \left\{
\begin{array}{ll} 
u_t = \nabla (c(|\nabla u|^2) \nabla u) & \mathrm{in} \ \Omega \times (0, + \infty) \\
\frac{\partial u}{\partial n} = 0 & \mathrm{in} \ \partial \Omega \times (0, + \infty) \\
u(x,0) = u_0 (x) & \mathrm{in} \ \Omega \\

\end{array} \right.
\end{equation}

Here $\Omega$ denotes picture domain. In general $\Omega$ is a bounded subset of $\mathbbm{R}^n$ with boundary of class $C^1$.  When looking for $C^1$ solutions, we assume that the derivatives of $u_0$ vanish at the boundary of $\Omega$. We usually restrict to diffusivity functions $c(s^2)$ monotonically decreasing from 1 to 0 while $s^2$ changes from 0 to $+ \infty$, infinitely smooth and such that function $\Phi(s) = sc(s^2)$ has one maximum in $\mathbbm{R}_{+}$ (although some authors investigated also different types of diffusivities). The typical choices are

\begin{equation}
\label{diffus}
 c(s^2) = \frac{1}{1+ \frac{s^2}{\lambda^2}} 
\end{equation}

and

\begin{equation}
\label{diffus2}
c(s^2) = \exp \left(- \frac{s^2}{2 \lambda^2}\right)
\end{equation}
with the parameter $\lambda > 0$. 
\\

\textit{The flux function} is defined as
\begin{equation}
 \Phi(s) = s \cdot c(s^2)
\end{equation}

For diffusivity \eqref{diffus} $\Phi '(s) < 0 $ for $|s| > \lambda$ and $\Phi '(s) > 0 $ for $|s| < \lambda$ (it is an easy calculation to give similar rule for 1.8).
\\

In the one-dimensional case \eqref{PMeq} simplifies to
\begin{equation}
 u_t  = \left( c(u^2_x) + 2u^2_x c'(u^2_x) \right) u_{xx}= \Phi '(u_x) u_{xx}
\end{equation}

\begin{figure}
\centering
\includegraphics[scale = 0.85]{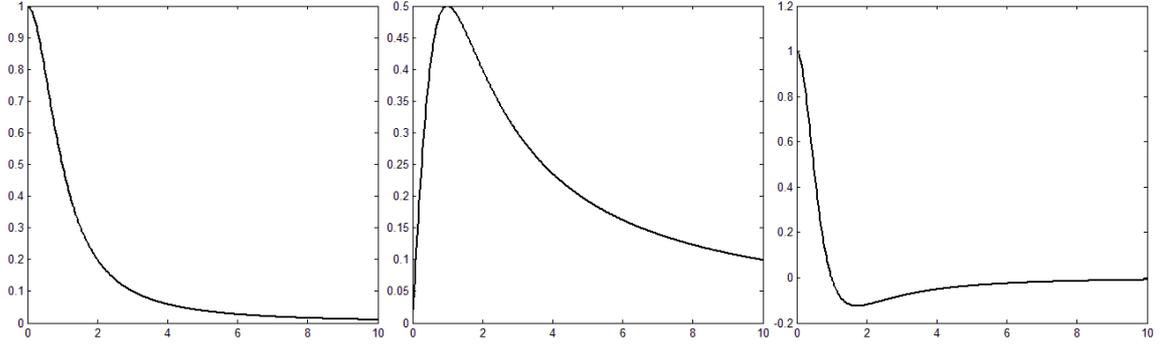}
\caption{Diffusivity function (1.7), the corresponding flux function and its derivative}
\end{figure}

Now it is visible, that for large values of $u_x$ the diffusion coefficient $\Phi '(u_x) $ becomes negative, which leads to backward diffusion.
\\

In the two-dimensional case \eqref{PMeq} can be rewritten as (see \cite{AGLM}):

\begin{equation}
\label{PMgauge}
u_t = \Phi ' (|\nabla u|) u_{\eta \eta} + c(|\nabla u|^2)u_{\xi \xi}
\end{equation}

where the coordinates $\eta$ and $\xi$ denote the directions parallel and perpendicular to $\nabla u$ respectively. Therefore diffusion is of forward type along the contour lines of the function $u$ -- \textit{isophotes} and forward-backward type in the perpendicular direction (gradient direction). We expect edges to be the regions of a large gradient values and therefore of backward diffusion in the gradient direction, which leads to sharpening instead of blurring. This explains how the Perona-Malik filter is not only able to prevent edges from being smoothed, but also enhance them.
\\

We see now that the Perona-Malik model is not anisotropic the in sense that diffusivity is a scalar function, not a tensor field. Nevertheless, since \eqref{PMgauge} shows that diffusivity coefficient varies in perpendicular directions, it is often referred to as an anisotropic diffusion equation.
\\
\\

\section{Investigation of well-posedness}
\ \\
 
\textbf{Definition 1.1} \textit{A \emph{well-posed} problem is a problem, for which there exists a unique solution, continuously depending on the data (in some reasonable topology). A problem, which is not well-posed is \emph{ill-posed}.}
\\

\textbf{Defintion 1.2} (Weak formulation in the one-dimensional case) \textit{ A locally integrable function $u(x,t)$ is said to be a \emph{weak solution} of the Perona-Malik equation if $\int_{\Omega} (u^2 + u_x^2) \, dx$ is uniformly bounded for bounded $t$ and if for any test function $\phi \in C^1_0(  \mathbbm{R}_+\times  \mathbbm{R})$}

\begin{equation}
\int \! \! \! \int [ \phi_t u - \phi_x c(u_x^2)u_{x} ] \, dx dt = 0
\end{equation}

Backward diffusion is well known to be an ill-posed process. Therefore possible presence of backward diffusion in Perona-Malik model suggests that the equation might be ill-posed. Pessimistic results came from H\"{o}llig \cite{Hol}, who gave an example of a different forward-backward diffusion process which can have infinitely many solutions. Since then many specific results were found. Below I survey some of the most important ones.
\\

Kichenassamy in \cite{Kit}, 1997,
showed that if any weak solution of the one-dimensional problem exists, the initial data  must be infinitely differentiable in the regions of backward diffusion ($|u_x| > \lambda $). This shows, that there might be no weak solution at all. Of course affine stationary functions ($u(x,t) = ax + b$) are solutions, but they are unstable meaning that there might be no solution if the initial condition is changed arbitrarily little. He also introduced the idea of looking for a generalized solutions of the Perona-Malik equation, showed that the Heavside step function is a generalized stationary solution and constructed a family of piecewise linear solutions.
\\

\textbf{Definition 1.3}\textit{ A function $u(x,t)$ is called a \emph{generalized solution} (or \emph{ultraweak solution}) of an equation $P(u) = 0$ in an open set D if there is a sequence $u_n$ of Lipschitz continuous functions such that $P(u_n) \rightarrow 0$ in the sense of distributions in D, and $u_n \rightarrow u$ in $L^1_{loc}(D)$.
A generalized solution of an evolution equation is called \emph{admissible} if its essential variation with respect to the space variables is nonincreasing in time.}
\\

\textbf{Theorem 1.2} (\cite{Kit} Properties of solutions)\textit{ There is no weak solution of problem \eqref{PMeq} in one dimension, that has bounded derivative which stays greater than $\lambda$ in a rectangle $(A,B) \times (0,T)$ for any $A,B \in \Omega$, $A < B$,  $T > 0$. }
\\

Kawohl and Kutev in \cite{KK}, 1998, proved that there are no global $C^1$ weak solutions of the one-dimensional problem with initial data involving backward diffusion (other than affine functions): 
\\

\textbf{Theorem 1.3} (\cite{KK} Nonexistance of solutions in the one-dimensional case) \textit{ Given $u_0$ such that $u'_0(x) = 0$ in $\partial \Omega$  and $u'_0(x) > \lambda$ (the diffusivity function parameter) in just one compact subinterval of $\Omega$ i.e. $u'_0(x) > \lambda$ in $ Q = (x_0,y_0) \subset \subset \Omega$ and $|u'_0(x)| < \lambda$ in $  \Omega / \overline{Q}$, $u_0$ analytic near $x_0$ (or strictly convex), the one-dimensional problem \eqref{PMeq} has no global weak solution in $C^1(\Omega)$. The result can be extended to finitely many subintervals and $u'_0(x) < - \lambda$.}
\\

They also established a useful maximum and comparison principles for the solutions and proved the uniqueness of local solutions (if they exist). Moreover, they showed that for initial data without backward diffusion there are classical global solutions in $C^2$.
\\

\textbf{Theorem 1.4} (\cite{KK} Maximum principle) \textit{ Suppose that $u$ is a Lipschitz continuous (weak) solution of \eqref{PMeq}. Then for every $p \in [2, \infty]$ the following inequality holds:}

\begin{equation}
||u(x,t)||_{L^p(\Omega)} \le ||u_0(x) ||_{L^p(\Omega)}
\end{equation}
\vspace{0.0mm}

\textbf{Theorem 1.5} (\cite{KK} Comparison principle, one-dimensional case)\textit{ Suppose that $u$ and $v$ are local weak solutions of \eqref{PMeq} with initial data $u_0$ and $v_0$. Then for $t \in (0,T): \ u_0(x) \le v_0(x) \Rightarrow u(x,t) \le v(x,t)$ if any of the following conditions holds:}
\\
a)
\[ \{ x \in \Omega : |u'_0(x)| > \lambda \} \cap \{ x \in \Omega : |v'_0(x)| > \lambda \} = \emptyset
\]
b)\begin{eqnarray*}
&& \exists w_0(x) \in C^{2,\alpha} \ \ \mbox{such that}\\
&& \forall x \in \Omega \ u_0(x) \le w_0(x) \le v_0(x), \ |w_0(x)| < \lambda, \ w'_0(x) = 0 \ \ \mbox{in} \ \ \partial \Omega \\
\end{eqnarray*}

\textbf{Theorem 1.6} (\cite{KK} Uniqueness of local weak solutions, one-dimensional case) \textit{ Suppose that $u$ and $v$ are local weak solutions of \eqref{PMeq} in $Q_T = \Omega \times (0,T)$ with identical analytic initial data $u_0$, $(u'_0)^2 - \lambda^2$ has only simple zeroes and diffusivity is an analytic function. Then $u(x,t) \equiv v(x,t)$ in $Q_T$.}
\\

More recently (2006) Rosati and Schiaffino showed \cite{RS} that piecewise linear interpolations of solutions of some discretizations of the Perona-Malik problem converge to ultraweak solutions.
\\

In 2009 Ghisi and Gobbino showed (for diffusivity \eqref{diffus}) that the set of initial data for  which there exists a local-in-time classical solution is dense in $C^1$ (\cite{GG}). They also proved an interesting fact that for dimensions higher than one Theorem 1.3 does not hold and the problem admits global solutions of the $C^{2,1}$ class even for the initial conditions indicating backward diffusion (\cite{GG2}).
\\

\section{Regularizations}

The idea behind any regularization of the Perona-Malik equation is to change the basic problem just slightly to obtain a well-posed problem. 
\\

In \cite{CLMC} Catte, Lions, Morel and Coll proposed using $u$ spatially convolved with the Gaussian function with standard deviation $\sigma$ in the diffusion coefficient. This leads to the following problem:

\begin{equation}
\label{regular}
 \left\{
\begin{array}{ll} 
u_t = \nabla (c(|\nabla G_{\sigma} * u|^2) \nabla u) = \nabla (c(|\nabla  u_{\sigma}|^2) \nabla u)  & \mathrm{in} \ \Omega \times (0, T] \\
\frac{\partial u}{\partial n} = 0 & \mathrm{in} \ \partial \Omega \times (0, T] \\
u(x,0) = u_0 (x) & \mathrm{in} \ \Omega \\
\end{array} \right.
\end{equation}

This problem turns into the Perona-Malik problem for $ \sigma= 0 $, if we assume $G_0(x) = \delta(x)$ (Dirac's delta distribution). For $\sigma > 0$ there is an unique regular solution of \eqref{regular} (a proof can be found in \cite{CLMC}).
\\

The intuition is that a solution exists because of the smoothing properties of convolution. Although this model became quite common for technical applications,  using the convolution with the Gaussian function is questionable, as Perona-Malik process was first proposed to replace Gaussian filter. 
\\

Different approaches involve regularization in time by averaging gradient value within some time interval, introducing a relaxation time into diffusivity \cite{Bar}, or mixing spatial and time regularization. 
\\

Interesting anisotropic regularization was given by Weickert [Weick], who proposed using a diffusion tensor $C$ instead of a scalar function. $C$ should have two eigenvectors

\[ v_1 \Vert  \nabla u_{\sigma} , \ \ \ \ v_2 \bot \nabla u_{\sigma} \]

with the eigenvalues 

\[ \lambda_1(\nabla u_{\sigma}) =  c(|\nabla  u_{\sigma}|^2), \ \ \ \ \lambda_2(\nabla u_{\sigma}) = 1 \]
\\
The eigenvalue acts as a diffusivity in the respective eigenvector direction.
Again, we obtain the Perona-Malik filter as $\sigma \rightarrow 0$ and $C \nabla u \rightarrow c(|\nabla u|^2) \nabla u $.

\chapter{Discretizations of the problem}\label{r:losers}

In this chapter I present the results obtained by Weickert \cite{Weick}. These explain how discretization turns the Perona-Malik equation into a well-posed problem.
\\

\textbf{Definition 2.1}\textit{ We call a problem originating from partial differential equation \emph{semidiscrete} if spatial variables are discrete and time remains continuous. A problem is \emph{discrete} (or \emph{fully-discrete}) if the time variable is also discrete.}

\section{Semidiscrete and discrete models}

One possible spatial discretization of \eqref{PMeq} in two dimensions can be written as (the discussion on discretizing parabolic equations can be found in \cite{MM}), i and j stand for vertical and horizontal pixel indices and with $\approx$ I denote approximating the derivatives by the difference quotients:

\begin{equation}
\label{discret}
 \left\{
\begin{array}{ll} 
c_{i,j} = c \left( \left( \frac{u_{i+1,j} - u_{i-1,j}}{2\Delta x} \right)^2 + \left( \frac{u_{i,j+1} - u_{i,j-1}}{2\Delta y} \right)^2 \right)\\
\\
c^{}_{i+1/2,j} := \frac{c^{}_{i,j} + c^{}_{i+1,j}}{2} \\
\\
\left[ c \frac{\partial u}{\partial x} \right]^{}_{i+1/2,j} \approx c^{}_{i+1/2,j} \left( \frac{u^{}_{i+1,j} - u^{}_{i,j}}{\Delta x} \right) := \varphi^{}_{i+1/2,j} \\
\\
\left[ \frac{\partial }{\partial x} \left( c \frac{\partial u}{\partial x} \right) \right]^{}_{i,j} \approx \frac{\varphi^{}_{i+1/2,j} - \varphi^{}_{i-1/2,j}}{\Delta x} \\
\\
\frac{d u_{i,j}}{d t} = \left[ \frac{\partial }{\partial x} \left( c \frac{\partial u}{\partial x} \right) \right]^{}_{i,j} + \left[ \frac{\partial }{\partial y} \left( c \frac{\partial u}{\partial y} \right) \right]^{}_{i,j} \approx \frac{\varphi^{}_{i+1/2,j} - \varphi^{}_{i-1/2,j}}{\Delta x} + \frac{\varphi^{}_{i,j+1/2} - \varphi^{}_{i,j-1/2}}{\Delta y}\\
\\
\frac{d u_{i,j}}{d t} \approx  \frac{1}{2 (\Delta y)^2} \left[(c_{i,j} + c_{i,j+1})(u_{i,j+1} - u_{i,j}) - (c_{i,j-1} + c_{i,j})(u_{i,j} - u_{i,j-1}) \right] + \ldots\\
\\
\ldots +\frac{1}{2 (\Delta x)^2} \left[(c_{i,j} + c_{i+1,j})(u_{i+1,j} - u_{i,j}) - (c_{i-1,j} + c_{i,j})(u_{i,j} - u_{i-1,j}) \right]\\

\end{array} \right.
\end{equation}
\\


Therefore, for the Neuman boundary condition, with $N(i,j)$ denoting the indices of elements in the neighbourhood of the element $u_{i,j}$ in the respective direction (a set of two elements or a single element), we have:
\begin{equation}
\frac{d u_{i,j}}{d t} = \sum_{(p,q) \in N_x(i,j)}\frac{c_{i,j} + c_{p,q}}{2(\Delta x)^2}(u_{p,q} - u_{i,j}) \ + \sum_{(p,q) \in N_y(i,j)} \frac{c_{i,j} +  c_{p,q}}{2 (\Delta y)^2} (u_{p,q} - u_{i,j})
\end{equation}

For $u_{i,j}$ in interior of the picture this equals:
\begin{equation}
\frac{d u_{i,j}}{d t} = \frac{1}{2 (\Delta x)^2} [ \overbrace{ (c_{i,j} + c_{i-1,j})}^{\beta_{i,j}^{-}} u_{i-1,j} \overbrace{- (c_{i-1,j} + 2c_{i,j} + c_{i+1,j})}^{\alpha_{i,j}} u_{i,j} + \overbrace{(c_{i,j} + c_{i+1,j})}^{\beta_{i,j}^{+}} u_{i+1,j} ] +\ldots
\end{equation}
\[ \ldots +  \frac{1}{2 (\Delta y)^2} [ \underbrace{(c_{i,j} + c_{i,j-1})}_{\delta_{i,j}^{-}} u_{i,j-1} \underbrace{-(c_{i,j-1} + 2c_{i,j} + c_{i,j+1})}_{\gamma_{i,j}} u_{i,j} + \underbrace{(c_{i,j} + c_{i,j+1})}_{\delta_{i,j}^{+}} u_{i,j+1} ]  \]
\\

Now, if we reorder $u$ to create a vector of length $M \cdot N$ (size of the image), we obtain:

\begin{equation}
\label{matrixA}
\frac{d \mathbf{u}}{dt} = (A_x + A_y) \mathbf{u} = A \mathbf{u}
\end{equation}

with $A$ - a sparse matrix of the size $MN \times MN$ with no more than five non-zero elements in each row ($\alpha_{ij} + \gamma_{ij},  \ \beta_{ij}^{+}, \ \beta_{ij}^{-}, \ \delta_{ij}^{+}, \  \delta_{ij}^{-}$).
\\

This semidiscrete model can be easily discretized in time, to obtain a fully discrete formulation, by putting

\begin{equation}
 \left\{
\begin{array}{ll} 
\frac{d u_{i,j}}{d t}(n \tau) \approx \frac{u^{n+1}_{ij} - u^{n}_{ij}}{\tau}\\
\\
u^{n+1}_{ij} = u^{n}_{ij} + \tau \left( \frac{\varphi^{n}_{i+1/2,j} - \varphi^{n}_{i-1/2,j}}{\Delta x} + \frac{\varphi^{n}_{i,j+1/2} - \varphi^{n}_{i,j-1/2}}{\Delta y} \right)\\

\end{array} \right.
\end{equation}

The scheme above is of explicit type, $\tau$ represents timestep. In a compact notation it can be rewritten as
\begin{equation}
\label{fulldiscret}
\mathbf{u}^{n+1} = \left[ I + \tau A(\mathbf{u}^n) \right] \mathbf{u}^n
\end{equation}
\\

In \cite{WRV} also a so-called semi-implicit scheme is analysed
\begin{equation}
\label{semiimp}
\mathbf{u}^{n} = \left[ I - \tau A(\mathbf{u}^n) \right] \mathbf{u}^{n+1}
\end{equation}
however in this case for each timestep we need to solve a system of equations, and the matrix $A$, although sparse, would have $2,56 \cdot 10^{10}$ elements for a picture of size $400 \times 400$ pixels.
\\

I will present here one more scheme, which is not a direct discretization of the Perona-Malik equation, but of a quite similar problem \eqref{PMani}. It is important, as it was proposed by Perona and Malik in the paper introducing their equation \cite{PM} and may serve as a simple example of an anisotropic filter.

\begin{equation}
\label{PMani}
u_t = \nabla \left( \left[ \begin{array}{rc}
c(u_x^2)&0\\
0&c(u_y^2)\\
\end{array} \right] \nabla u \right)
\end{equation}

The difference between the problems is that the Perona-Malik equation would have the same expression, $c(|\nabla u|^2)$, on the whole diagonal. 
\\
\vspace{15 mm}

The discretization (explicit scheme) proposed by Perona and Malik is

\begin{equation}
 \left\{
\begin{array}{ll} 
c_{i+1/2,j}^n := c \left( \left( \frac{u_{i+1,j}^n - u_{i,j}^n }{\Delta x} \right)^2 \right)\\
\\
c_{i,j+1/2}^n := c \left( \left( \frac{u_{i,j+1}^n - u_{i,j}^n }{\Delta y} \right)^2 \right)\\
\\
\frac{ u_{i,j}^{n+1} - u_{i,j}^{n}}{ \tau} =  \frac{1}{2 (\Delta x)^2} \left[ c_{i+1/2,j}^{n}(u_{i+1,j}^n - u_{i,j}^n) - c_{i-1/2,j}^{n}(u_{i,j}^n - u_{i-1,j}^n) \right] + \ldots \\
\\
 \ldots + \frac{1}{2 (\Delta y)^2} \left[ c_{i,j+1/2}^{n}(u_{i,j+1}^n - u_{i,j}^n) - c_{i,j-1/2}^{n}(u_{i,j}^n - u_{i,j-1}^n) \right] \\

\end{array} \right.
\end{equation}



\section{Well-posedness of discrete models}

The discretization \eqref{discret} leads to a nonlinear system of ordinary differential equations

\begin{equation}
 \left\{
\begin{array}{ll} 
\frac{d \mathbf{u}}{dt} = A(\mathbf{u}) \mathbf{u}\\
\\
\mathbf{u}(0) = \mathbf{u_0} \\
\end{array} \right.
\end{equation}

In \cite{Weick} a general analysis of such semidiscrete problems is given. Here I present only the main results in the form of a theorem.
\\

\textbf{Theorem 2.1} (\cite{Weick} Properties of semidiscrete problems) \textit{If the matrix $A$ of a semidiscrete problem in the form \eqref{discret} possesses properties P1 - P5, then the described process satisfies the properties I1 - I4.}
\\

\textit{
\textbf{P1} $A(\mathbf{u})$ is a continuously differentiable function of $\mathbf{u}$.}
\\

\textit{
\textbf{P2}  $A(\mathbf{u})$ is a symmetric matrix.}
\\

\textit{
\textbf{P3}  All row sums of $A(\mathbf{u})$ are zero.}
\\

\textit{
\textbf{P4}  All off-diagonal elements of $A(\mathbf{u})$ are nonnegative.}
\\

\textit{
\textbf{P5}  $A(\mathbf{u})$ is irreducible, which means that for any two indexes $i, j$ there exist a sequence of indices $(k_0,k_1, \dots , k_r)$ such that $k_0 = i, k_r = j, a_{k_p,k_{p+1}} \neq 0$ for $p = 0,1, \ldots, r-1$.}
\\

Matrix $A(\mathbf{u})$ obtained in \eqref{matrixA} possesses the above properties, which is an easy observation.
\\

\textit{
\textbf{I1} The problem is well-posed, meaning there is unique solution $\mathbf{u}(t) \in C1([0,T], \mathbbm{R}^{MN})$ for every $T > 0$, depending continuously on the initial value and the right-hand side of
the ODE system (stability).} This is implied by P1 (with the Picard-Lindel\"{o}f theorem).
\\


\textit{
\textbf{I2}  Average grey level invariance (this is a usefull property in practical applications)}

\[ \forall t > 0 : \ \ \frac{1}{MN} \sum_{k =1}^{MN} (\mathbf{u_{0}})_k = \frac{1}{MN} \sum_{k =1}^{MN} (\mathbf{u}(t))_k = \mu \]

This is implied by P2 and P3.
\\

\textit{
\textbf{I3} The extremum principle}
\[ \forall t > 0, \ \ \forall i \in J: \ \ \min_{k \in J} (\mathbf{u_0})_k  \le (\mathbf{u}(t))_i \le  \max_{k \in J} (\mathbf{u_0})_k  , \ \ J = \{1,2, \ldots , MN \} \]
This is implied by P3 and P4.
\\

\textit{
\textbf{I4}  Convergence to a constant steady state}
\[ \forall i \in J: \ \ \lim_{t \rightarrow \infty} (\mathbf{u}(t))_i = \mu \]
This is implied by P2-P5, moreover, certain functions (like the variance of data) are decreasing througout the process.
\\

Analysis of the fully discrete case is also given in \cite{Weick}. Once again, I present only the most important conclusion.
\\

\textbf{Theorem 2.2} (\cite{Weick} Properties of discrete problems) \textit{ Scheme \eqref{fulldiscret} inherits the properties of the semidiscrete scheme \eqref{discret} (I1 - I4) for} 
\[ \tau < \frac{1}{\frac{2}{(\Delta x)^2} + \frac{2}{(\Delta y)^2}} \]

\textit{
which is a stability restriction. Scheme \eqref{semiimp} is unconditionally stable and therefore always possesses the properties I1-I4. }
\\

Theorem 2.2 shows that the discretization of the ill-posed Perona-Malik equation may lead to a~well-posed problem and therefore explains the Perona-Malik paradox by showing the stability of numerical implementations.

\chapter{Numerical results}

\section{Comparison between Perona-Malik and Gaussian filters}

Figure 3.1 presents the evolution of one dimensional data under Gaussian (left column) and Perona-Malik (right column) filters - explicit scheme \eqref{discret}. 

\begin{figure}[H]
\centering
\includegraphics[scale = 0.75]{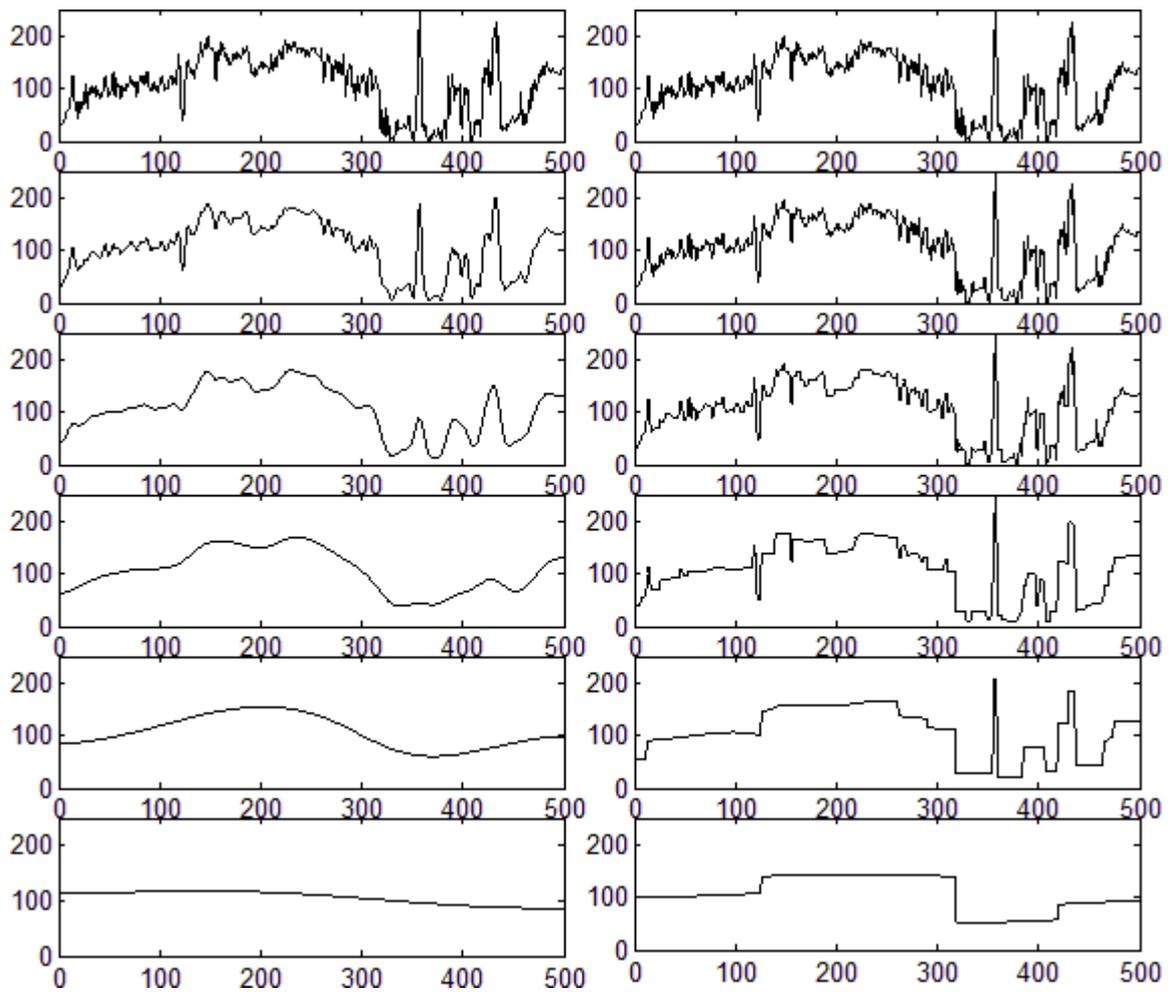}
\caption{Initial data and evolution for 10, 100 , 1000, 10000, 100000 iterations with the timestep 0.1 and $ \lambda = 1$.}
\end{figure}

It is visible on 3.1 how the Perona-Malik filter prevents structures like edges from being smoothed. Another observation is the tendency to develop small regions of constant values of the function.

\section{Different diffusivity functions}

The following figure presents the differences in the evolution of data for different diffusivity functions: model \eqref{diffus} (on the left) and \eqref{diffus2} (on the right) for the explicit scheme \eqref{discret}.

\begin{figure}[H]
\centering
\includegraphics[scale = 0.8]{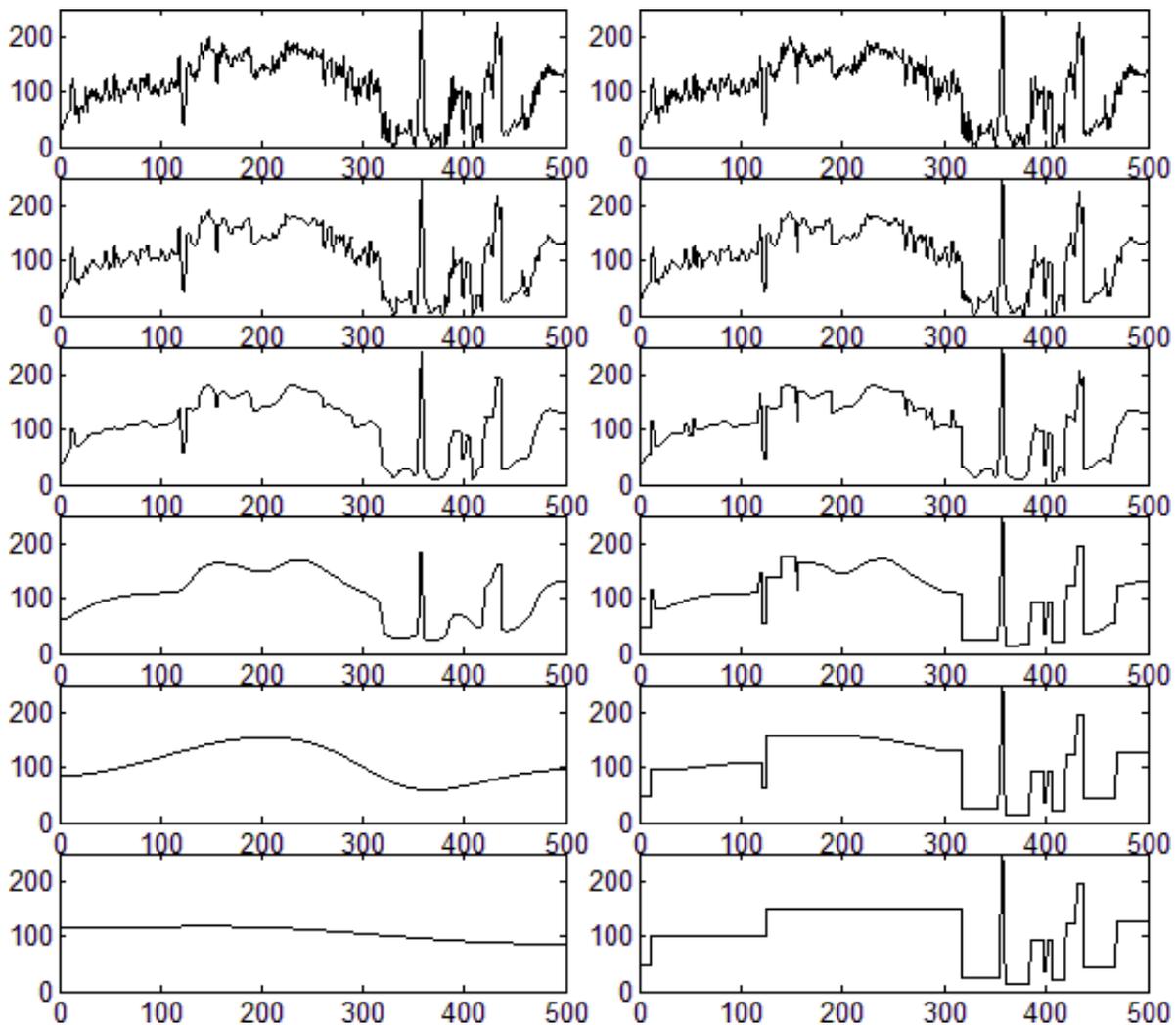}
\caption{Initial data and evolution for 10, 100 , 1000, 10000, 100000 iterations with timestep 0.1 and $ \lambda = 5$.}
\end{figure}
There is a significant difference. The first model acts more like a smoothing filter and flattens the data faster, while the second one preserves the structure for a longer time. This is due to larger negative values of the flux derivative and the enhancing properties of backward diffusion process (fig. 3.3). 
\\
\\

\begin{figure}[H]
\centering
\includegraphics[scale = 0.7]{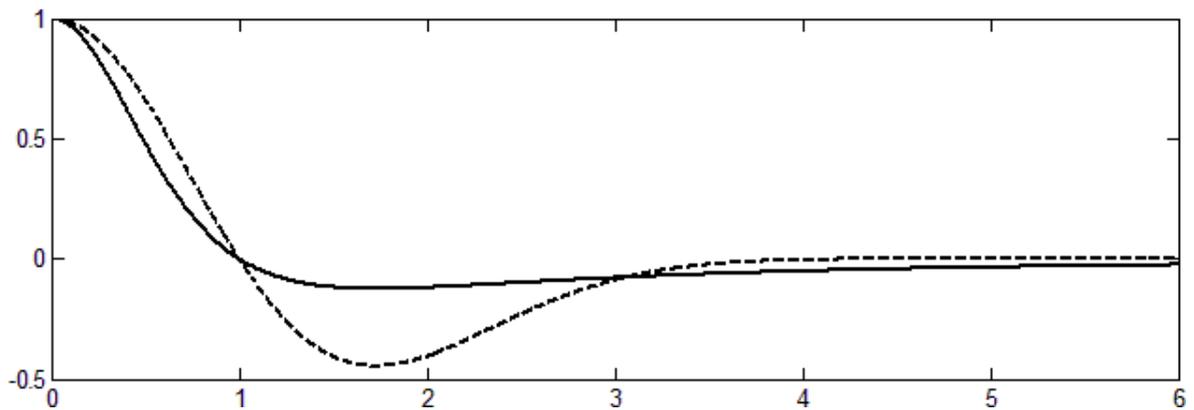}
\caption{Derivatives of the flux function for diffusivities \eqref{diffus} (continuous) and \eqref{diffus2} (dotted) for $\lambda = 1$}
\end{figure}

\begin{figure}[H]
\centering
\includegraphics[scale = 0.65]{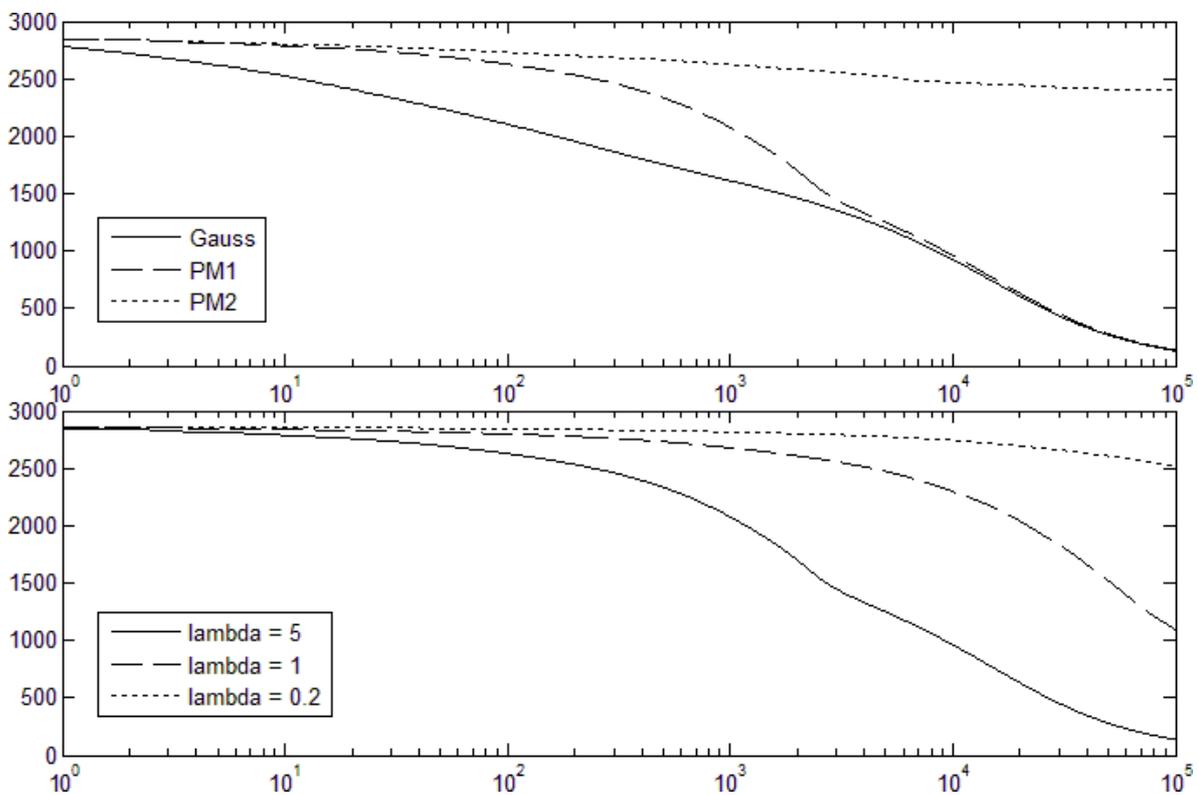}
\caption{Data variance as the function of the number of iterations for processes on fig. 3.2 (PM1 and PM2, respectively) and Gaussian filtration and (lower plot) for diffusivity \eqref{diffus} with different $\lambda$}
\end{figure}

\section{2D data filtration and evolution}

In this section I compare the Perona-Malik filter with the anisotropic process \eqref{PMani} and the Perona-Malik filter regularized as in \eqref{regular} with $\sigma = 1$ (pixel). For testing the filtration capability I use a picture with additional Gaussian noise and measure the distance from the real picture using the metric
\begin{equation}
d(u_0,u) = \sum |u_0(i,j) - u(i,j)|
\end{equation}

Although filtration is not the basic application of the Perona-Malik filter, it is an application for which results are easily measurable. The following figures present the relative error and variance as the function of the number of iterations. Initial data had SNR ratio equal 2 (signal to noise ratio, measured by the mean signal divided by the noise standard deviation).  All algorithms were stopped in the first minimum of the error function. The parameter $\lambda$ was arbitrary chosen to be equal 1, timestep 0.2, diffusivity function \eqref{diffus}. It took 12690 iterations for the aniostropic algorithm to reach the minimum, 3833 iterations for the standard algorithm and 1557 for the regularized one. Clearly, the regularized algorithm needed much less iterations than the other two and gave slightly better results.

\begin{figure}[H]
\centering
\includegraphics[scale = 0.6]{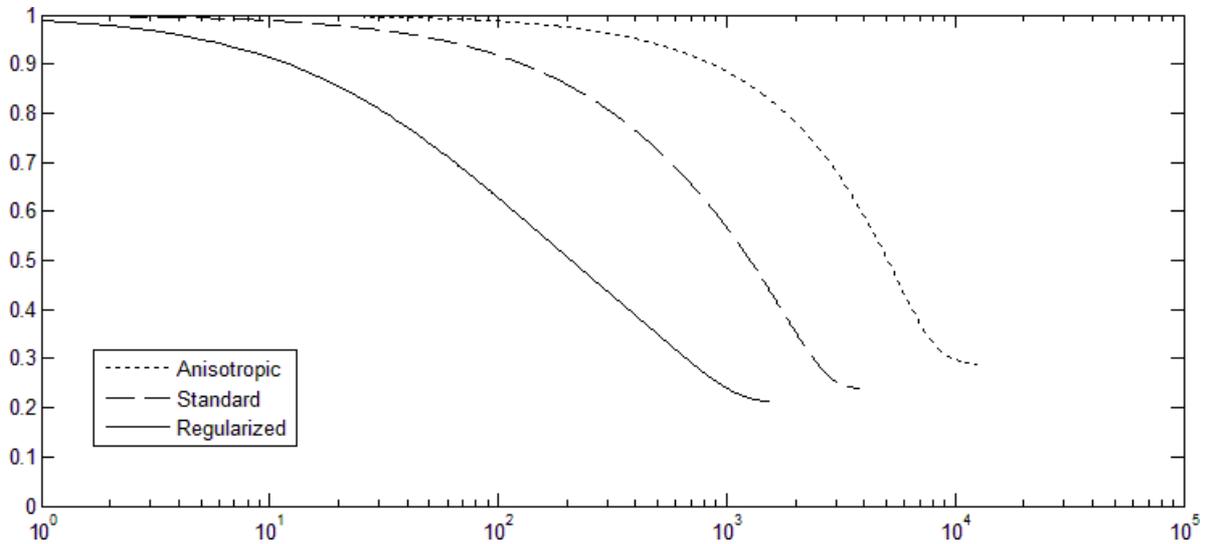}
\caption{Relative error as the function of the number of iterations}
\end{figure}

\begin{figure}[H]
\centering
\includegraphics[scale = 0.59]{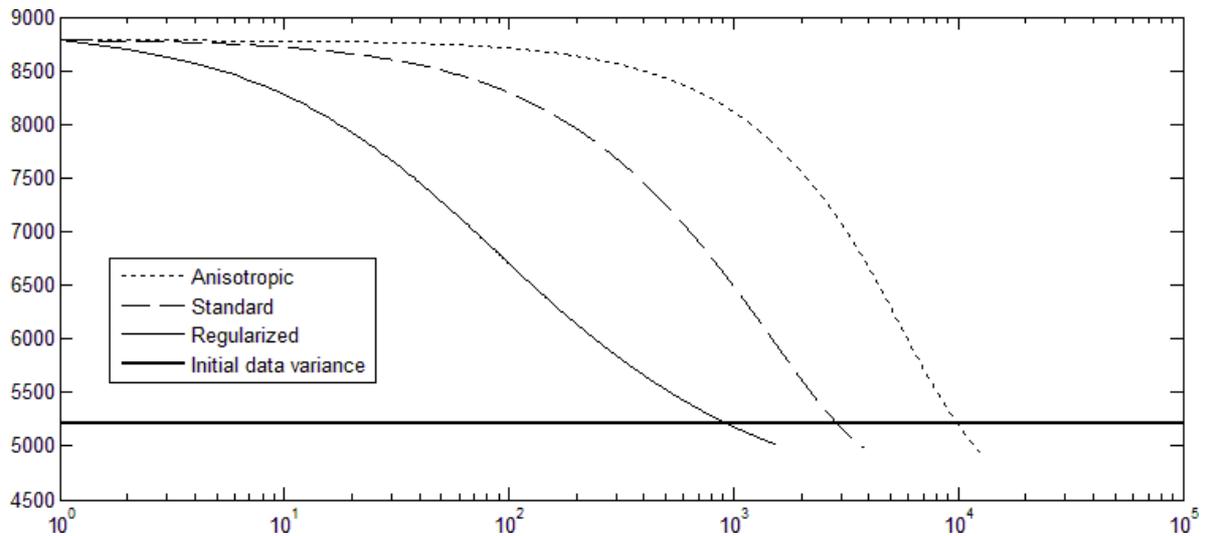}
\caption{Variance as the function of the number of iterations}
\end{figure}

\begin{figure}[H]
\centering
\includegraphics[scale = 1]{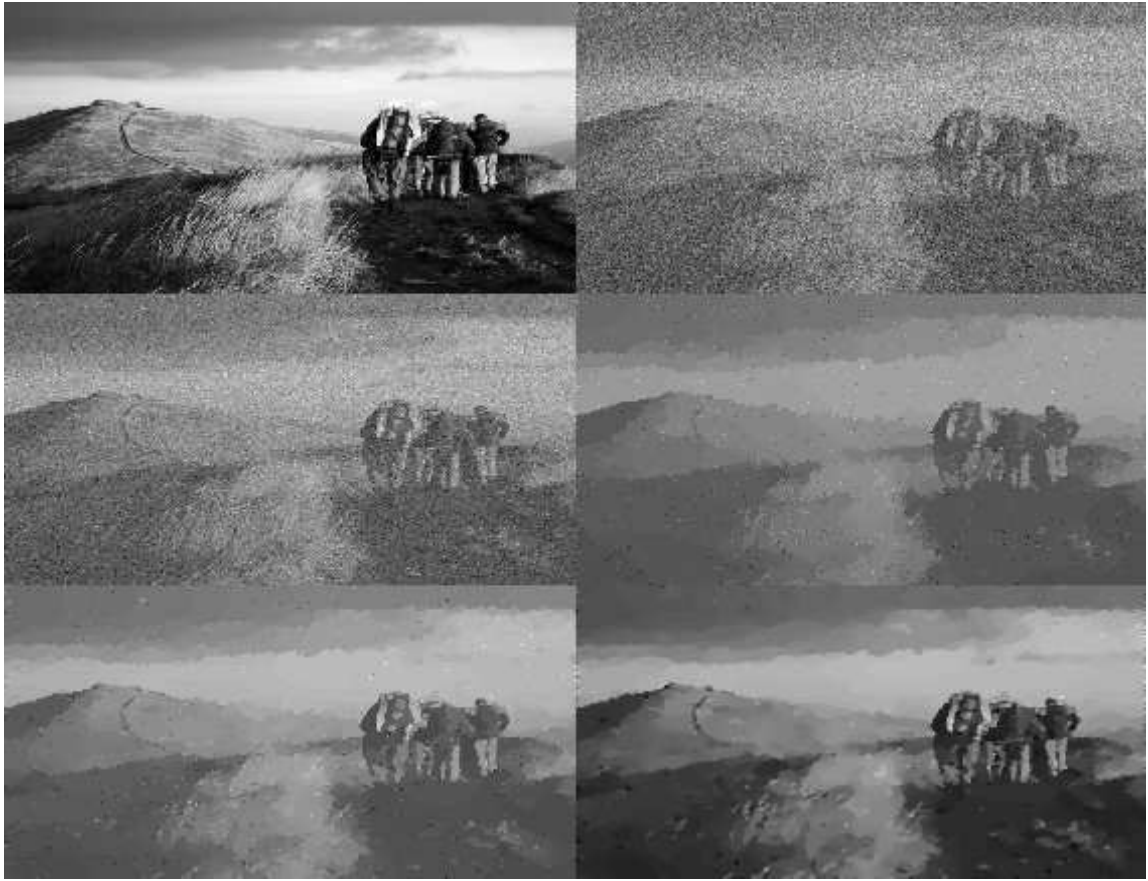}
\caption{Left to right: picture, picture with noise, picture after 1000 iterations with the standard algorithm, result with the anisotropic algorithm, result with the standard algorithm, result with the regularized algorithm}
\end{figure}

\begin{figure}[H]
\centering
\includegraphics[scale = 0.59]{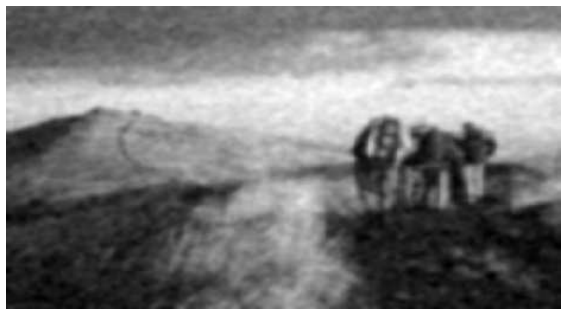}
\caption{For comparison, the effect of Gaussian filtration}
\end{figure}

\begin{figure}[H]
\centering
\includegraphics[scale = 1.0]{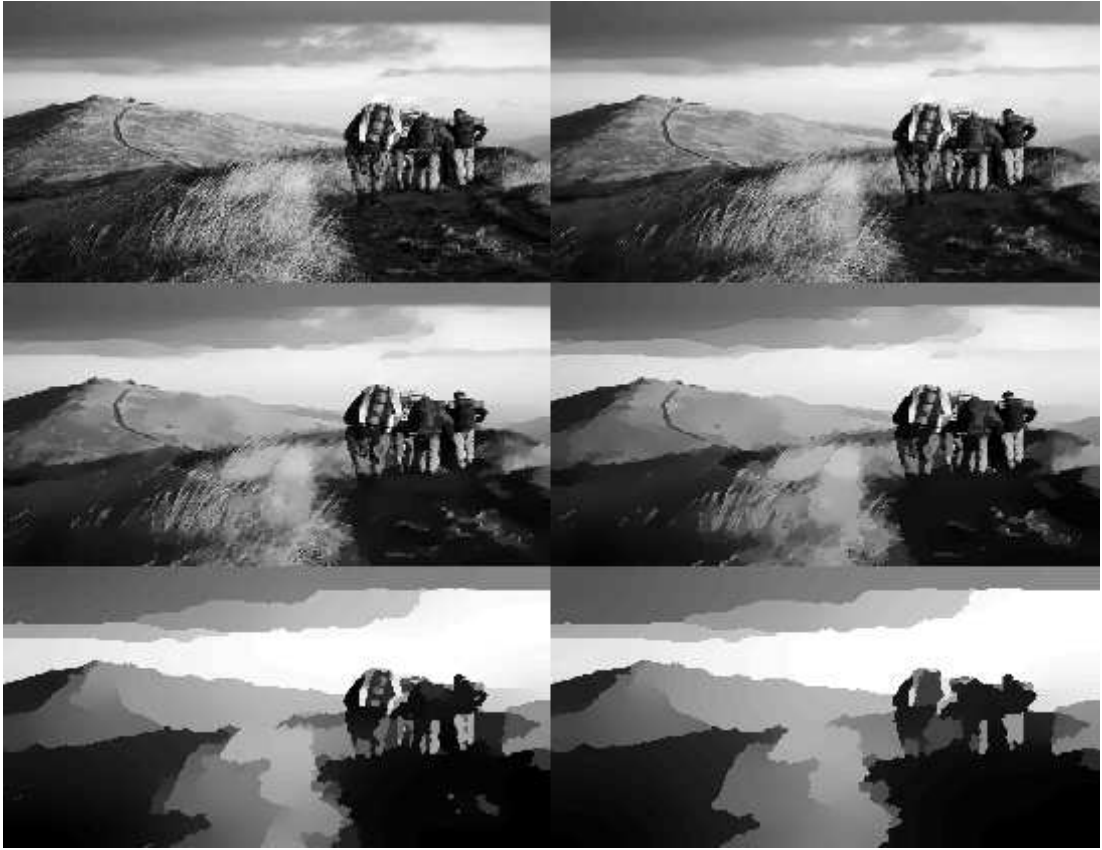}
\caption{Evolution of the image "Polonina" under the Perona-Malik filter for 0, 100, 500, 1000, 5000 and 10000 iterations with timestep 0.2 and $\lambda = 1$}
\end{figure}

It is visible, how the Perona-Malik filter may serve as a tool for image segmentation.
\\

In this chapter as the two-dimensional data I used image "Polonina", bitmap created from a photograph taken by me in Bieszczady mountains in 2009 and as the one-dimensional data I used one horizontal line taken from the image "Polonina".

\appendix

\chapter{List of implemented diffusion algorithms}

\ \\ 

\textbf{PMexp} - Perona-Malik filter with explicit scheme for two-dimensional data
\\

\textbf{PMexp1D} - Perona-Malik filter with explicit scheme for one-dimensional data
\\

\textbf{PMexpR} - regularized Perona-Malik filter with explicit scheme for two-dimensional data
\\

\textbf{PeronaMalik} - Perona-Malik filter in anisotropic version, 2D
\\

\textbf{PMimp} - Perona-Malik filter with semi-implicit scheme, 2D
\\

\textbf{PMimp1D} - Perona-Malik filter with semi-implicit scheme, 1D
\\

\textbf{GaussDiffusion} - simple gaussian diffusion
\\


\begin{thebibliography}{99}
\addcontentsline{toc}{chapter}{Bibliography}

\bibitem[AGLM]{AGLM} Luis Alvarez, Frederic Guichard, Pierre-Louis Lions, Jean-Michel Morel, \textit{ Axioms and fundamental equations in
image processing}, Archive for Rational Mechanics and Analysis, vol. 123 (1993), 199-257.


\bibitem[Bar]{Bar} G. I. Barenblatt, M. Bertsch, R. Dal Passo, M. Ughi, \textit{ A degenerate pseudoparabolic regularization of a nonlinear forward-backward heat equation arising in the theory of heat and mass exchange in stably stratified turbulent shear flow}, SIAM Journal on Mathematical Analysis, vol. 24 (1993),
1414-1439.

\bibitem[CLMC]{CLMC} Francine Catté, Pierre-Louis Lions, Jean-Michel Morel, Tomeu Coll,   \textit{Image Selective Smoothing and Edge Detection by Nonlinear Diffusion}, SIAM Journal on Numerical Analysis, vol. 29 no. 1. (1992), 182-193.

\bibitem[GG]{GG} Marina Ghisi, Massimo Gobbino, \textit{A class of local classical solutions for the one-dimensional Perona-Malik equation}, Journal: Transaction of the American Mathematical Society 361 (2009), 6429-6446.

\bibitem[GG2]{GG2} Marina Ghisi, Massimo Gobbino, \textit{An example of global classical solution for the Perona-Malik equation}, arXiv:0907.0772v1 (2009).

\bibitem[Hol]{Hol} Klaus Hollig,
  \textit{Existence of Infinitely Many Solutions for a Forward Backward Heat Equation}, Transactions of the American Mathematical Society, , vol. 278. no. 1. (1983), 299-316.

\bibitem[Kit]{Kit} Satyanad Kitchenassamy, \textit{The Perona-Malik Paradox}, SIAM Journal on Applied Mathematics, vol. 57 no. 5. (1997), 1328-1342.

\bibitem[KK]{KK}Bernd Kawohl, Nikolai Kutev, \textit{Maximum and comparison principles for one-dimensional anisotropic diffusion}, Mathematische
Annalen 311. (1998), 107-123.

\bibitem[MM]{MM} K. W. Morton, D. F. Mayers, \textit{Numerical Solution of Partial Differential Equations}, Cambridge University Press (1994).

\bibitem[PM]{PM} Pietro Perona, Jitendra Malik, \textit{Scale-Space and Edge Detection Using Anisotropic
Diffusion}, IEEE Transactions on Pattern Analysis and Machine Intelligence, vol. 12. no. 7. (1990), 629-639.


\bibitem[RS]{RS} Mario Rosatia, Andrea Schiaffinob, \textit{Some remarks about Perona-Malik equation}, Nonlinear Analysis 65 (2006), 1-11.

\bibitem[Strzel]{Strzel} Paweł Strzelecki, \textit{Krótkie wprowadzenie do równań różniczkowych cz\k{a}stkowych}, Wydawnictwa Uniwersytetu Warszawskiego, Warsaw, Poland (2006).

\bibitem[WB]{WB} Joachim Weickert, Brahim Benhamouda, \textit{A Semidiscrete Nonlinear Scale-Space Theory and Its Relation to the Perona-Malik Paradox}, Theoretical Foundations of Computer Vision (1997), 1-10.

\bibitem[Weick]{Weick} Joachim Weickert, \textit{Anisotropic diffusion in image processing}, Teubner-Verlag, Stuttgart, Germany (1998).

\bibitem[WRV]{WRV} Joachim Weickert, Bart M. ter Haar Romeny, Max A. Viergever, \textit{Efficient and Reliable Schemes for Nonlinear Diffusion Filtering}, IEEE Transactions on Image Processing, vol. 7. no. 3. (1998), 398-410.

\end{thebibliography}
\end{document}